\documentclass[10pt,a4paper]{amsart}
\usepackage[utf8]{inputenc}
\usepackage[T1]{fontenc}
\usepackage{%
	lmodern,%
	mparhack,%
	relsize%
} 
\usepackage[english]{babel}
\usepackage{%
	amsaddr,%
	amsfonts,%
	amsmath,%
	amsthm,%
	amscd,%
	amssymb,%
	mathtools%
}
\usepackage{xcolor}   
\usepackage[hyphens]{url}
\usepackage{hyperref}
\hypersetup{
	breaklinks=true,%
	colorlinks=true,
	frenchlinks=true,%
	linkcolor=red,%
	citecolor=red,%
	filecolor=magenta,%
	urlcolor=magenta,%
	linktoc=all,
}
\usepackage{orcidlink} 	      
\usepackage[final]{microtype} 
\usepackage[%
	backend=biber,%
    citestyle=numeric-comp,%
	style=numeric-comp,%
	sorting=nty,%
	hyperref,%
	block=space,%
	backref%
]{biblatex}
\usepackage[autostyle=true]{csquotes}

\newcommand{\exampleref}[1]{\hyperref[example:#1]{example~\ref{example:#1}}}
\newcommand{\algref}[1]{\hyperref[alg:#1]{algorithm~\ref{alg:#1}}}
\newcommand{\chapref}[1]{\hyperref[chap:#1]{chapter~\ref{chap:#1}}}
\newcommand{\secref}[1]{\hyperref[sec:#1]{section~\ref{sec:#1}}}
\newcommand{\subsecref}[1]{\hyperref[subsec:#1]{subsec.~\ref{subsec:#1}}}
\newcommand{\subsubsecref}[1]{\hyperref[subsubsec:#1]{subsubsec.~\ref{subsubsec:#1}}}
\newcommand{\figref}[1]{\hyperref[fig:#1]{fig.~\ref{fig:#1}}}
\newcommand{\tabref}[1]{\hyperref[tab:#1]{table~\ref{tab:#1}}}
\newcommand{\defref}[1]{\hyperref[def:#1]{definition~\ref{def:#1}}}
\newcommand{\thmref}[1]{\hyperref[thm:#1]{theorem~\ref{thm:#1}}}
\newcommand{\lemmaref}[1]{\hyperref[lemma:#1]{lemma~\ref{lemma:#1}}}
\newcommand{\corollaryref}[1]{\hyperref[corollary:#1]{corollary~\ref{corollary:#1}}}
\newcommand{\remarkref}[1]{\hyperref[rmk:#1]{remark~\ref{rmk:#1}}}

\makeatletter
\newcommand{\homcoords}[1]{\left[#1\checknextarg}
\newcommand{\checknextarg}{\@ifnextchar\bgroup{\gobblenextarg}{\right]}}
\newcommand{\gobblenextarg}[1]{ : #1\@ifnextchar\bgroup{\gobblenextarg}{\right]}}
\makeatother

\newcommand{\commandof}[2]{\mathmbox{#1\roundbrack{#2}}}
\newcommand{\keyphrase}[1]{\textcolor{blue}{#1}\index{#1}}

\newcommand{\roundbrack}[1]{\mathmbox{\left(#1\right)}}
\newcommand{\modp}[1]{\ \mathmbox{\roundbrack{\mathrm{mod}\ #1}}}
\newcommand{\divp}[1]{\mathmbox{\mathrm{div}\roundbrack{#1}}}
\newcommand{\Ker}{\mathrm{Ker}}
\newcommand{\Kerof}[1]{\commandof{\Ker}{#1}}
\newcommand{\Img}{\mathrm{Im}}
\newcommand{\Imgof}[1]{\commandof{\Img}{#1}}
\newcommand{\Mod}{\mathrm{Mod}}

\newcommand{\Modx}[1]{\Mod_{#1}}
\newcommand{\Modxof}[2]{\commandof{\Modx{#1}}{#2}}

\newcommand{\Exp}{\mathrm{Exp}}

\newcommand{\Expx}[1]{{\Exp}_{#1}}
\newcommand{\Expxof}[2]{\commandof{\Expx{#1}}{#2}}
\newcommand{\Expxinv}[1]{{\Exp}^{-1}_{#1}}
\newcommand{\Expxinvof}[2]{\commandof{\Expxinv{#1}}{#2}}

\newcommand{\cardof}[1]{\commandof{\mathrm{card}}{#1}}
\newcommand{\C}{\mathbb{C}}

\newcommand{\Q}{\mathbb{Q}}
\newcommand{\Z}{\mathbb{Z}}

\newcommand{\K}{\mathbb{K}}
\newcommand{\Zof}[1]{\frac{\Z}{#1 \Z}}
\newcommand{\weier}{\mathcal{W}}
\newcommand{\redweier}{\overline{\mathcal{W}}}
\newcommand{\edwards}{\mathcal{E}}

\newcommand{\jacobian}{\mathcal{J}}
\newcommand{\jacobianof}[1]{\commandof{\jacobian}{#1}}
\newcommand{\qpjacobianof}[2]{\commandof{\jacobian_{#1}}{#2}}
\newcommand{\qpjacobianzof}[2]{\commandof{\jacobian_{#1}^{0}}{#2}}
\newcommand{\jacobianzof}[1]{\commandof{\jacobian^{0}}{#1}}

\newcommand{\fullcurveof}[2]{\commandof{#1}{#2}}
\newcommand{\weierof}[1]{\fullcurveof{\weier}{#1}}
\newcommand{\redweierof}[1]{\fullcurveof{\redweier}{#1}}
\newcommand{\edwardsof}[1]{\fullcurveof{\edwards}{#1}}
\newcommand{\Qp}[1]{\Q_{#1}}
\newcommand{\GF}[1]{\commandof{\mathrm{GF}}{#1}}

\newtheorem{theorem}{Theorem}
\numberwithin{theorem}{section} 

\newtheorem{remark}{Remark}
\numberwithin{remark}{section} 
\newtheorem{corollary}[theorem]{Corollary}
\newtheorem{definition}{Definition}
\numberwithin{definition}{section} 
\newtheorem{example}{Example}

\numberwithin{equation}{section}

\title{Exp function for Edwards curves over local fields}
\author{Giuseppe Filippone \orcidlink{0000-0001-7315-1852}}
\address{%
	\normalsize Department of Mathematics and Computer Science\\%
	\normalsize University of Palermo\\%
	\normalsize Via Archirafi 34, 90123 Palermo, Italy}
\email{\href{mailto:giuseppe.filippone01@unipa.it}{giuseppe.filippone01@unipa.it}}
\date{}

\begin{document}

	\begin{abstract}
		We extend the exponential map $ \Exp $ for complex elliptic curves in short Weierstrass form to
		Edwards curves over local fields. Subsequently, we compute the map $ \Exp $ for Edwards curves
		over the local field of $ p $-adic numbers.
	\end{abstract}

	\maketitle

	\section{Introduction}\label{introduction}
	The literature on elliptic curves and their applications in \keyphrase{cryptography} is well consolidated.
	Recently, curves such as Montgomery elliptic curves and Edwards curves
	(in particular in their \keyphrase{twisted} version) have gained great popularity
	for their cryptographic applications.

	Edwards curves were first introduced in $ 2007 $ by H. Edwards \cite{Edwards}.
	These curves are already the subject of many papers in
	cryptography
	\cite{Lange2011, BernsteinLange, BernsteinLange2, HisilWongCarterDawson, BernsteinBirknerLangePeters}.
	Compared to the classic elliptic curves in Weierstrass form, they are more efficient
	for cryptographic use and the (single or multiple) digital signature.
	An application of Edwards curves to Goppa Codes is shown in \cite{GFGoppaCodes}.

	Since the Weierstrass elliptic functions fulfill the identity
	$ {\roundbrack{\frac{1}{2} \wp^\prime(z)}}^2 = \wp^3(z) - \frac{g_2}{4} \wp(z) - \frac{g_3}{4} $,
	where $ g_2, g_3 \in \C $ are constants, the function
	$ \Exp \colon z \mapsto \roundbrack{\wp(z), \frac{1}{2} \wp^\prime(z)} $
	maps an element $ z $ belonging to the complex torus $ \C / \Lambda $,
	where $ \Lambda $ is the period lattice of $ \wp $, to a point belonging to
	the corresponding elliptic curve in short Weierstrass form of the complex projective plane,
	defined by the equation $ y^2 = x^3 - \frac{g_2}{4} x - \frac{g_3}{4} $.
	Moreover, it is such that
	$ \Exp(z_1 + z_2) = \Exp(z_1) \ast \Exp(z_2) $ (see e.g. $ \S $VI and $ \S $IX in \cite{Silverman2009}),
	where the operation $ \ast $
	is given by the chord-and-tangent law on the points of the elliptic curve.

	In this paper, we extend the above exponential map to
	\keyphrase{Edwards curves} over \keyphrase{local fields},
	and we give a particular specialization of this map over
	the local field $ \Qp{p} $ of $ p $-adic numbers.
	We are motivated by authoritative literature on the matter of lifting,
	summarized in \cite{SilvermanLiftingECDLP}
	where the author gives a survey connecting the lifting to the
	\keyphrase{discrete logarithm problem} over elliptic curves in Weierstrass form.

	Although cryptosystems over infinite fields have received little attention in the past,
	in \cite{Xu2008} the authors gave a cryptosystem based on
	quotient groups of an elliptic curve in Weierstrass form over the
	$ p $-adic number field, able to encrypt messages
	with variable lengths. This led to public-key cryptosystems
	with hierarchy management \cite{XuYue},
	which look interesting for their possible applications.

	More recently, similar topics have been investigated in \cite{Tang2015},
	where the authors consider twisted Edwards curves over local fields and introduce a cryptosystem
	based on quotient groups of twisted Edwards curves over local fields.

	For these reasons, although it is possible to extend the above map to other forms of
	elliptic curves (such as Legendre form, Jacobi form, Hessian form,
	Huff form), in this work we will focus only on
	the Edwards form.
	To the best of our knowledge,
	there are no other papers in which this study was already addressed.

	In \secref{1}, we describe Edwards curves and their relationship with elliptic curves in Weierstrass form.
	In \secref{2}, we extend the map $ \Exp $ for elliptic curves in Weierstrass form over $ \C $ to
	Edwards curves over local fields.
	Finally, in \secref{3}, we exhibit the map $ \Exp $ for Edwards curves when the
	local field taken into account is the field $ \Qp{p} $ of $ p $-adic numbers.

	\section{Prerequisites and notations}\label{sec:1}

	The goal of this paper is to compute the map $ \Exp $ for the Edwards curves over local fields.
	For a general introduction to local fields, we address the reader to a classic book,
	e.g. \cite{SerreLocalField}.
	Here we summarize some results on Edwards curves, which will be used later
	and give explicitly a reduction (\thmref{jacobian-edwards}) to canonical forms
	of divisors on an Edwards curve and an explicit equivalence, under particular conditions, between
	a class of Edwards curves and a class of elliptic curves in Weierstrass form (see
	\thmref{isomorphism-w-e}).

	\begin{definition}[Edwards curves]
		A (non-smooth) algebraic curve over a field $ \K $ which,
		with respect to a suitable coordinate system, has the equation
		$ {\hat{x}}^2 + {\hat{y}}^2 = 1 + d {\hat{x}}^2 {\hat{y}}^2 $,
		where $ d \in \K $ is such that $ d (d - 1) \ne 0 $,
		is called an Edwards curve $ \edwards $.
	\end{definition}

	Recall that, over a field $ \K $ of characteristic different from $ 2 $, a (smooth) elliptic curve
	(possessing at least a $ \K $-rational point) can be represented in a suitable coordinate
	system by the Weierstrass equation $ y^2 = x^3 + a^\prime x^2 + b^\prime x $, having one point
	at infinity $ \Omega = [Z : X : Y] = [0: 0: 1] $ on the $ y $-axis. Hence, from here on, unless
	otherwise specified, we will consider an elliptic curve in Weierstrass form defined by the latter
	equation.

	\begin{remark}
		Note that, unlike those in Weierstrass form, curves in Edwards form $ \edwards $
		have two points at infinity, that is,
		$ \Omega_1 = \mathmbox{[\hat{Z} : \hat{X} : \hat{Y}]} = \mathmbox{[0: 1: 0]} $ on
		the $ x $-axis and $ \Omega_2 = \mathmbox{[\hat{Z} : \hat{X} : \hat{Y}]} = \mathmbox{[0: 0: 1]} $
		on the $ y $-axis, which are \keyphrase{ordinary singular points} for
		$ \edwards $.
	\end{remark}

	In the following, we provide a brief introduction to the group law for Edwards
	curves, which was first considered in \cite{BernsteinLange} (cf. also \cite{Edwards,ARENE2011842}).

	Formally speaking, one has to take into account the group of divisor classes
	$ \frac{\operatorname{Div}^0(\edwards)}{\operatorname{Princ}(\edwards)} $,
	modulo the subgroup	of principal divisors on $ \edwards $. In particular,
	one wants to sum the two divisors $ (P - O) $ and $ (Q - O) $,
	where $ P, Q \in \edwardsof{\K} $ are two affine points of $ \edwards $, and
	$ O = (0, 1) \in \edwardsof{\K} $ is taken as the base point, in light
	of \thmref{jacobian-edwards}.

	Let $ \kappa $ be the unique hyperbola (containing $ O^\prime = (0, -1), 2 \Omega_1 $ and $ 2 \Omega_2 $),
	passing through $ P $ and $ Q $, which intersects the curve $ \edwards $ in a further point $ R = (x_R, y_R) $.
	Let $ l_R \colon Y - y_R Z = 0 $ be the line passing through $ R $ and parallel to the $ x $-axis (thus $ l_R $
	pass through $ S = (-x_R, y_R) $ as well).
	One has that
	\begin{equation*}
		\mathrm{div}\left(\frac{\kappa}{(Y - y_R Z) \cdot X}\right) = P + Q - S - O,
	\end{equation*}
	hence $ \mathmbox{(P - O) + (Q - O) \equiv (S - O)} $.

	The above group law can be summarized into the following
	addition and doubling formulas, where for all (not necessarily distinct) points
	$ P = (x_P, y_P) $ and $ Q = (x_Q, y_Q) $, the sum divisor $ S - O \equiv (P - O) + (Q - O) $ is
	such that:
	\begin{equation*}
		S = \roundbrack{\frac{x_P y_Q + x_Q y_P}{1 + d x_P x_Q y_P y_Q}, \frac{y_P y_Q - x_P x_Q}{1 - d x_P x_Q y_P y_Q}}.
	\end{equation*}

	\begin{remark}[cf. {\normalfont \cite{BernsteinLange}}]\label{rmk:opposite}
		Note that $ (P - O) + (Q - O) \equiv O - O $ if and only if
		$ P = (x_P, y_P) $ and $ Q = (x_Q, y_Q) $ are such that $ y_P = y_Q $ and $ x_P = -x_Q $, that is,
		$ Q $ is the symmetric point, with respect to the $ y $-axis, to the point $ P $.
	\end{remark}

	\begin{remark}\label{rmk:d-non-square}
		Note that if the parameter $ d $ is a non-square, then the denominators in the addition
		and doubling formulas cannot vanish {\normalfont \cite{BernsteinLange}} and the affine points
		of the curve give in turn a subgroup of the whole group of divisor classes
		(cf. \corollaryref{subgroup-affine}).
	\end{remark}

	In terms of the group of divisor classes, one finds either of the following reduced divisors
	in each divisor class.

	\begin{theorem}[Jacobian of Edwards curves]\label{thm:jacobian-edwards}
		Let $ \edwards $ be an Edwards curve, and let $ \jacobianof{\edwards} $ be the Jacobian
		of $ \edwards $.
		Every divisor $ D \in \jacobianof{\edwards} $ has one of the following \keyphrase{canonical forms}:
		\begin{enumerate}
			\itemsep0em
			\item $ D \equiv P - O $;
			\item $ D \equiv (P - O) + (\Omega_1 - O) $;
			\item $ D \equiv (P - O) + (\Omega_2 - O) $;
			\item $ D \equiv (P - O) + (\Omega_1 - \Omega_2) $,
		\end{enumerate}
		where $ P \in \edwards $ is an affine point.
		In particular, the divisors equivalent to $ P - O $ form a subgroup $ \jacobianzof{\edwards} $
		(\corollaryref{subgroup-affine}) of index $ 4 $ in $ \jacobianof{\edwards} $, and
		$ 2 \Omega_1 \equiv O^\prime + O $ and $ 2 \Omega_2 \equiv H^\prime + H $, where
		$ O = (0, 1) $, $ O^\prime = (0, -1) $, $ H = (1, 0) $, and $ H^\prime = (-1, 0) $.
	\end{theorem}
	\begin{proof}
		Let $ D = D_1 + D_2 $ be a divisor of $ \edwards $, where $ D_1 $
		is such that every point in its support is an affine point,
		and $ D_2 = t_1 \Omega_1 + t_2 \Omega_2 $ with $ t_1, t_2 \in \Z $.

		We show that every even multiple of $ \Omega_1 $ and $ \Omega_2 $ is
		equivalent to a multiple of $ O^\prime + O $ and $ H + H^\prime $, respectively.
		Indeed, we have that:
		\begin{equation*}
			\begin{aligned}
				\divp{\frac{X}{Z}} &= O^\prime + O - 2 \Omega_1,\\
				\divp{\frac{Y}{Z}} &= H^\prime + H - 2 \Omega_2,\\
			\end{aligned}
		\end{equation*}
		thus $ O^\prime + O \equiv 2 \Omega_1 $ and $ H^\prime + H \equiv 2 \Omega_2 $.

		As a consequence, we can reduce $ D $ to one of the canonical forms shown in the claim,
		by exploiting the above rule, the group law for $ \edwards $ and the following remark:
		if $ t_1 $ and $ t_2 $ are both odd, first we reduce $ D $ to one of these two forms
		$ (P - O) + (\Omega_1 - \Omega_2) $ or $ (P - O) + (\Omega_2 - \Omega_1) $, because
		$ D $ is a zero degree divisor, but the latter is equivalent to the former because
		$ \Omega_1 - \Omega_2 \equiv (\Omega_2 - \Omega_1) + (O^\prime + O) - (H^\prime + H) $.

		Finally, since
		\begin{equation*}
			\begin{aligned}
				2 (\Omega_1 - O) &= 2\Omega_1 - 2 O \equiv (O^\prime + O) - 2 O = O^\prime - O \in \jacobianzof{\edwards},\\
				2 (\Omega_2 - O) &= 2\Omega_2 - 2 O \equiv H^\prime + H - 2 O \equiv O - O \in \jacobianzof{\edwards},
			\end{aligned}
		\end{equation*}
		the quotient group $ \frac{\jacobianof{\edwards}}{\jacobianzof{\edwards}} $ is isomorphic
		to $ \frac{\Z}{2 \Z} \oplus \frac{\Z}{2 \Z} $.
	\end{proof}

	\begin{corollary}\label{corollary:subgroup-affine}
		The subset $ \jacobianzof{\edwards} $ of zero degree divisors whose support contains only affine points is a subgroup
		of $ \jacobianof{\edwards} $.
	\end{corollary}
	\begin{proof}
		Let $ (P - O), (Q - O) \in \jacobianzof{\edwards} $ be two divisors, where $ P, Q \in \edwardsof{\K} $.
		By \remarkref{opposite} one has that $ -(Q - O) \in \jacobianzof{\edwards} $,
		and by \remarkref{d-non-square} one has that $ (P - O) - (Q - O) \in \jacobianzof{\edwards} $.
	\end{proof}

	\begin{remark}
		Given an elliptic curve in Weierstrass form $ \weier $,
		it is usual to identify the non-zero divisor $ P - \Omega $ with the point
		$ P $ of $ \weierof{\K} $, and the zero divisor $ \Omega - \Omega $ with $ \Omega $.
		Similarly, one can denote the non-zero divisor $ P - O $ of $ \jacobianzof{\edwards} $
		with the affine point $ P $ of $ \edwardsof{\K} $, and the zero divisor
		$ O - O $ with $ O $.
		Hence, one may refer to either the Jacobian or the group of
		$ \K $-rational points of these curves, indifferently.
	\end{remark}

	In the following, we describe under which conditions one has an equivalence between Edwards curves
	$ \edwards $ and elliptic curves in Weierstrass form $ \weier $.

	\begin{definition}[cf. \cite{BernsteinLange, GFGoppaCodes}]\label{def:w-e-bi-eq}
		Let $ \edwards $ be an Edwards curve defined,
		over a field $ \K $ of characteristic different from $ 2 $, by the equation
		$ {\hat{x}}^2 + {\hat{y}}^2 = 1 + d {\hat{x}}^2 {\hat{y}}^2 $, where $ d (d - 1) \ne 0 $.
		Let $ 0 \ne x_1 \in \K $ be such that $ x_1 $ and $ (1 - d) $
		are both non-square or square in $ \K $, and let $ y_1 \in \K $ such that
		$ y_1^2 = \frac{4 x_1^3}{1 - d} $.
		Putting $ a^\prime = 2 x_1 \frac{1 + d}{1 - d} $ and $ b^\prime = x_1^2 $,
		one considers the elliptic curve in Weierstrass form $ \weier = \weier_{d, x_1} $, defined over $ \K $
		by the equation	$ y^2 = x^3 + a^\prime x^2 + b^\prime x $, and one denotes by $ \alpha $
		and $ \beta $ the two following rational maps:
		\begin{subequations}
			\begin{align}
				\alpha\colon \edwardsof{\K} &\longrightarrow \weierof{\K} \label{eq:alpha} \\ \nonumber
				(\hat{x}, \hat{y}) &\longmapsto (x, y) = \left(x_1 \frac{1 + \hat{y}}{1 - \hat{y}}, y_1 \frac{(1 + \hat{y})}{\hat{x} (1 - \hat{y})}\right),\\
				\alpha^{-1} = \beta \colon \weierof{\K} &\longrightarrow \edwardsof{\K} \label{eq:beta} \\ \nonumber
				(x, y) &\longmapsto (\hat{x}, \hat{y}) = \left(\frac{y_1 x}{x_1 y}, \frac{x - x_1}{x + x_1}\right),
			\end{align}
		\end{subequations}
		which make $ \weier $ and $ \edwards $ birationally equivalent.

		Moreover, one extends the definition of $ \alpha $ and $ \beta $  by putting $ \alpha((0, 1)) = \Omega $,
		$ \beta(\Omega) = (0, 1) $, $ \alpha((0, -1)) = (0, 0) $ and $ \beta((0, 0)) = (0, -1) $; and
		(possibly) $ \beta((t_1, 0)) = \beta((t_2, 0)) = \Omega_1 $, $ \beta((-x_1, \pm s_1)) = \Omega_2 $,
		where $ (t_1, 0), (t_2, 0), (-x_1, \pm s_1) \in \weierof{\K} $, with $ t_1, t_2 \ne 0 $ (see
		{\normalfont \cite{GFGoppaCodes}} for further comments).
	\end{definition}

	In the following, we stress the meaning of taking $ d $ a non-square in a field $ \K $.

	\begin{theorem}[Isomorphism between $ \jacobianof{\weier} $	and $ \jacobianzof{\edwards} $]\label{thm:isomorphism-w-e}
		Let both $ \edwards $ and $ \weier $ as in \defref{w-e-bi-eq}. If $ d $ is
		not a square in the field $ \K $, then
		there is an isomorphism over $ \K $ between the group $ \jacobianof{\weier} $
		and the subgroup $ \jacobianzof{\edwards} $ defined in \thmref{jacobian-edwards}.
	\end{theorem}

	\begin{proof}
		By using the parameters $ x_1 $, $ y_1 $, $ a^\prime $, and $ b^\prime $
		as defined in \defref{w-e-bi-eq}, we prove that, if $ d $ is a non-square,
		then the rational map $ \beta $ in \eqref{eq:beta} defines a biregular map
		between the elliptic curve in Weierstrass form $ \weier $ of equation
		$ y^2 = x^3 + a^\prime x^2 + b^\prime x $ and the subset of $ \edwards $
		consisting of its affine points. We are left with proving that there is no point in
		$ \weier $ with abscissa $ -x_1 $ and that $ (0, 0) $ is the only point in $ \weierof{\K} $
		with ordinate $ y = 0 $.

		The former assertion follows from the fact that, by intersecting
		the line $ x = - x_1 $ and the curve $ \weier $ one has that
		$ -x_1^3 + a^\prime x_1^2 - b^\prime x_1 $ is equal to $ d y_1^2 $,
		which is a non-square in $ \K $ because $ d $ is a non-square and, therefore, there is
		no point in $ \weierof{\K} $ with abscissa $ - x_1 $.

	 	The latter assertion follows from the fact that
		the intersection between the line $ y = 0 $ and the curve $ \weier $
		has no roots in $ \K $ except $ x = 0 $. More precisely, since
		$ \Delta(x^2 + a^\prime x + b^\prime) = d {\roundbrack{\frac{4 x_1}{1 - d}}}^2 $
		is not a square in $ \K $ because $ d $ is a non-square, then
		$ (0, 0) $ is the only point in $ \weierof{\K} $
		with ordinate $ y = 0 $.

		Since the map $ \beta $ in {\normalfont \eqref{eq:beta}} transforms
		a line through $ P \in \weierof{\K} $ and $ Q \in \weierof{\K} $
		onto the hyperbola through $ \beta(P) $, $ \beta(Q) $,
		$ O^\prime, 2 \Omega_1 $ and $ 2\Omega_2 $, and maps vertical lines onto horizontal lines,
		then $ \beta $ induces a group homomorphism of the corresponding divisor classes groups.
	\end{proof}

	From here on, we confine ourselves to the case in \thmref{isomorphism-w-e}, where one can simply identify
	the elements $ P - O $ of $ \jacobianzof{\edwards} $ with the point $ P $, so that
	$ O $ is the neutral element of the group.

	\section{The map \texorpdfstring{$ \Exp $}{Exp} for Edwards curves over local fields}\label{sec:2}

	In this section, we use the already known results on elliptic curves in Weierstrass form,
	and we extend these results to Edwards curves.

	Recall that, as we said in the previous section,
	we have confined ourselves to the case where there is a birational equivalence
	between an elliptic curve $ \weier $ in Weierstrass form and an Edwards curve $ \edwards $
	such that the Jacobian $ \jacobianof{\weier} $ of $ \weier $ is isomorphic to
	$ \jacobianzof{\edwards} $ (\thmref{isomorphism-w-e}), that is,
	the subgroup of divisors of $ \edwards $ whose reduced
	form is $ P - O $, where $ P $ is an affine point and $ O = (0, 1) $ is taken as
	the neutral element of the group.

	Let $ \K $ be a local field, $ \mathcal{O}_\K $ its ring of integers,
	$ \mathfrak{m}_\K $ its prime ideal, and $ k = \mathcal{O}_\K / \mathfrak{m}_\K $
	its residue field.
	We will take the image $ \qpjacobianzof{k}{\edwards} $ under reduction modulo $ \mathfrak{m}_\K $ of the
	group $ \qpjacobianzof{\K}{\edwards} $, then we will investigate under what assumptions
	one has that
	$ \qpjacobianzof{\K}{\edwards} \cong \qpjacobianzof{k}{\edwards} \oplus \mathfrak{m}_\K $.

	First, we remark that for elliptic curves $ \redweier $ in short Weierstrass form, defined
	by the equation $ y^2 = x^3 + a x + b $, whose
	reduction modulo $ \mathfrak{m}_\K $ is non-singular, the following sequence:
	\begin{equation}
		0 \longrightarrow \mathfrak{m}_\K \xrightarrow{\Expx{\redweier}} \qpjacobianof{\K}{\redweier} \ \xrightarrow{\Modx{\redweier}} \qpjacobianof{k}{\redweier} \longrightarrow 0
	\end{equation}
	is exact \cite{kosterspannekoek} (see also \cite[ch. $ \S $VII]{Silverman2009}),
	thus $ \Imgof{\Expx{\redweier}} = \Kerof{\Modx{\redweier}} $,
	$ \Expx{\redweier} $ is a monomorphism, $ \Modx{\redweier} $ is an epimorphism, and one has that
	\begin{equation*}
		\qpjacobianof{k}{\redweier} \cong \frac{\qpjacobianof{\K}{\redweier}}{\Kerof{\Modx{\redweier}}} = \frac{\qpjacobianof{\K}{\redweier}}{\Imgof{\Expx{\redweier}}}.
	\end{equation*}

	The map $ \Modx{\redweier} $ is nothing else than a simple reduction modulo $ \mathfrak{m}_\K $ of the coordinates
	of the points $ P = \homcoords{Z}{X}{Y} $ in $ \redweierof{\K} $ which, up to a multiplication times a suitable
	$ t \in \mathcal{O}_\K $, have integral entries $ Z, X, Y \in \mathcal{O}_\K $:
	\begin{alignat*}{2}
		\Modx{\redweier} \colon \   &\qpjacobianof{\K}{\redweier}&\ \rightarrow\  &\qpjacobianof{k}{\redweier}\\
		&P = \homcoords{Z}{X}{Y} &\ \mapsto\  &\homcoords{Z \modp{\mathfrak{m}_\K}}{X \modp{\mathfrak{m}_\K}}{Y \modp{\mathfrak{m}_\K}}.
	\end{alignat*}
	Note that $ \Modx{\redweier} $ is trivially surjective for Hensel's lemma
	(see proof in \cite[sec. $ \S $VII.2.1]{Silverman2009}).

	Furthermore, the function $ \Expx{\redweier} $ is defined as follows:
	\begin{alignat*}{2}
		\Expx{\redweier} \colon \  &\mathfrak{m}_\K &\ \longrightarrow\  &\qpjacobianof{\K}{\redweier}\\
		&z &\ \longmapsto\  & \homcoords{1}{\wp(z)}{\frac{1}{2} {\wp^\prime(z)}}.
	\end{alignat*}

	\begin{remark}\label{rmk:series-expansion}
		Note that the Weierstrass $ \wp $-function (and its derivative) can be expressed {\normalfont \cite{AbramowitzStegun1972, ApostolTomMike1990MfaD}}
		through a Laurent series in a neighborhood of zero, and one has that
		\begin{align*}
			\wp(z) &= \frac{1}{z^2} + \sum_{k=2}^{\infty} c_k z^{2k - 2},\\
			\wp^\prime(z) &= -\frac{2}{z^3} + \sum_{k=2}^{\infty} (2k - 2) c_k z^{2k - 3},
		\end{align*}
		where $ c_2 = \frac{g_2}{20} $, $ c_4 = \frac{g_3}{28} $,
		$ c_k = \frac{3}{(2k + 1)(k - 3)} \sum_{m=2}^{k - 2} c_m c_{k - m} $,
		and $ g_2, g_3 \in \C $ are the parameters of the elliptic curve over $ \C $ in short Weierstrass
		form defined by the equation
		$ {\roundbrack{\frac{\wp^\prime(z)}{2}}}^2 = \wp^3(z) - \frac{g_2}{4} \wp(z) - \frac{g_3}{4} $.

		Moreover, one may generalize these results over a local field $ \K $ taking into account a
		neighborhood of zero in order to have a convergence for the series expansion of
		$ \wp $ and $ \wp^\prime $.
		In this case, one has that $ g_2 $, $ g_3 $ and $ z $ belong to the local field $ \K $.
	\end{remark}

	Since $ z = 0 $ is the only element of $ \mathfrak{m}_\K $
	mapped to $ \Omega $, the homomorphism $ \Expx{\redweier} $
	is into; thus, for any $ z $ in a neighborhood of zero, one can define
	\begin{equation}\label{eq:inv-exp}
		\Expxinv{\redweier} \coloneqq -2\frac{\wp(z)}{\wp^\prime(z)}
	\end{equation}
	(see $ \S $IV and $ \S $VII \cite{Silverman2009} for further details)
	such that $ \Expxinvof{\redweier}{\Expxof{\redweier}{z}} = z $, whose first terms in
	a Taylor series are
	\begin{equation*}
		z + \frac{g_2}{10} z^5 + \frac{3 g_3}{28} z^7 + \frac{g_2^2}{120} z^9 +
			\frac{23 g_2 g_3}{1540} z^{11} + O(z^{13}).
	\end{equation*}

	Thus, taking $ \edwards $ such that its reduction modulo $ \mathfrak{m}_\K $
	is non-singular, we have the following theorem.

	\begin{theorem}[The map $ \Exp $ for Edwards curves]\label{thm:exp-edwards}
		Let $ \K $ be a local field, $ \mathcal{O}_\K $ is its ring of integers,
		and $ \mathfrak{m}_\K $ is the prime ideal of $ \mathcal{O}_\K $.
		If $ \edwards $ is an Edwards curve as in \thmref{isomorphism-w-e}, that is, with
		$ d \in \K $ a non-square, then the following map:
		\begin{alignat*}{2}
			\Expx{\edwards} \colon\  &\mathfrak{m}_\K &\ \longrightarrow\  &\qpjacobianzof{\K}{\edwards}\\
				&z 		    &\ \longmapsto\  & \roundbrack{\frac{2}{3}\frac{y_1 (3 \wp(z) - a^\prime)}{x_1 \wp^\prime(z)}, \frac{3 \wp(z) - a^\prime - 3 x_1}{3 \wp(z) - a^\prime + 3 x_1}},
		\end{alignat*}
		where $ x_1 $, $ y_1 $, and $ a^\prime $ are as in \defref{w-e-bi-eq},
		is an exponential map for $ \edwards $, that is,
		$ \Expxof{\edwards}{z_1 + z_2} = \Expxof{\edwards}{z_1} + \Expxof{\edwards}{z_2} $.
	\end{theorem}

	\begin{proof}
		Recall that, in \defref{w-e-bi-eq}, we have a birational equivalence between the Edwards curve $ \edwards $
		and the elliptic curve in Weierstrass form $ \weier $ of
		equation $ y^2 = x^3 + a^\prime x^2 + b^\prime x $, whereas
		the above map $ \Expx{\redweier} $ is defined for elliptic curves
		in short Weierstrass form $ \redweier $ of equation $ y^2 = x^3 + ax + b $.

		However,  we can apply the transformation
		$ \chi \colon (x, y) \mapsto \roundbrack{x - \frac{a^\prime}{3}, y} $ which,
		through the change of variables $ \bar{x} = x - \frac{a^\prime}{3} $, $ \bar{y} = y $,
		changes the Weierstrass form $ y^2 = x^3 + a^\prime x^2 + b^\prime x $
		onto the short Weierstrass form $ {\bar{y}}^2 = {\bar{x}}^3 + a \bar{x} + b $, that is,
		for any $ P = (\bar{x}, \bar{y}) \in \redweierof{\K} $
		such that $ {\bar{y}}^2 = {\bar{x}}^3 + a \bar{x} + b $,
		we have that $ \chi(P) = P^\prime \in \weierof{\K} $.
		As $ \chi(P) = P^\prime $ belongs to $ \weierof{\K} $, 	we
		can now compute $ \beta(P^\prime) $, where $ \beta $ in \eqref{eq:beta}, in order to get a
		point belonging to $ \edwardsof{\K} $. In particular, if $ P = \Expxof{\redweier}{z} $ for
		some $ z \in \mathfrak{m}_\K $, then
		\begin{equation*}
			\begin{aligned}
				\beta\roundbrack{\chi\roundbrack{P}} =
				\beta\roundbrack{\chi\roundbrack{\Expxof{\redweier}{z}}} &=
				\beta\roundbrack{\chi\roundbrack{\homcoords{1}{\wp(z)}{\frac{1}{2} \wp^\prime(z)}}} = \\
				&= \beta\roundbrack{\homcoords{1}{\wp(z) - \frac{a^\prime}{3}}{\frac{1}{2} \wp^\prime(z)}} = \\
				&= \roundbrack{\frac{2}{3}\frac{y_1 (3 \wp(z) - a^\prime)}{x_1 \wp^\prime(z)}, \frac{3 \wp(z) - a^\prime - 3 x_1}{3 \wp(z) - a^\prime + 3 x_1}}.\\
			\end{aligned}
		\end{equation*}

		Thus, the map $ \Expx{\edwards} $ for Edwards curves over the local field $ \K $ is defined as
		$ \Expx{\edwards} = \beta \circ \chi \circ \Expx{\redweier} $, that is,
		\begin{alignat*}{2}
			\Expx{\edwards} \colon &\mathfrak{m}_\K \ &\longrightarrow\ &\qpjacobianof{\K}{\edwards}\\
			&z \ &\longmapsto\ &
				\Expxof{\edwards}{z}
				= \roundbrack{\frac{2}{3}\frac{y_1 (3 \wp(z) - a^\prime)}{x_1 \wp^\prime(z)}, \frac{3 \wp(z) - a^\prime - 3 x_1}{3 \wp(z) - a^\prime + 3 x_1}}.\\
		\end{alignat*}

		Note that $ \chi\roundbrack{\Omega} = \Omega = \homcoords{0}{0}{1} $ as the projective map
		$ \chi $ maps $ \homcoords{Z}{X}{Y} $ onto the point $ \homcoords{Z}{X - \frac{a^\prime}{3}Z}{Y} $, and
		thus we have that $ \beta\roundbrack{\chi\roundbrack{\Expxof{\redweier}{0}}} =
		\beta\roundbrack{\chi\roundbrack{\Omega}} = \beta(\Omega) = O $.

		Finally, we are left to prove that the map
		$ \Expx{\edwards} = \beta \circ \chi \circ \Expx{\redweier} $ is a one-to-one homomorphism of groups.
		On the one hand, the maps $ \Expx{\redweier} $, $ \beta $, and $ \chi $ are one-to-one.
		Indeed, the map $ \beta $ here is bijective as $ d $ is not a square,
		and $ \chi^{-1} \colon (\bar{x}, \bar{y}) \mapsto \roundbrack{x + \frac{a^\prime}{3}, y} $.

		On the other hand, $ \Expx{\edwards} $ is a homomorphism because $ \Expx{\redweier} $
		and $ \beta $ (see \thmref{isomorphism-w-e}) are homomorphisms, and
		$ \chi $ is a translation, thus one has that:
		\begin{equation*}
			\begin{aligned}
				\Expxof{\edwards}{z_1 + z_2} &= \beta \circ \chi \circ \Expxof{\redweier}{z_1 + z_2} = \\
				&= \beta \circ \chi \circ \roundbrack{\Expxof{\redweier}{z_1} + \Expxof{\redweier}{z_2}} = \\
				&= \beta \circ \chi \circ \Expxof{\redweier}{z_1} + \beta \circ \chi \circ \Expxof{\redweier}{z_2} =\\
				&= \Expxof{\edwards}{z_1} + \Expxof{\edwards}{z_2}.
			\end{aligned}
		\end{equation*}
	\end{proof}

	\begin{remark}
		As $ \chi $ transforms the curve $ \weier $ into the curve $ \redweier $,
		one has that $ \chi\roundbrack{P_1 + P_2} = \chi\roundbrack{P_1} + \chi\roundbrack{P_2} $,
		where the left term uses the addition formula for $ \weier $, and the right
		term uses the addition formula for $ \redweier $.
	\end{remark}

	\begin{corollary}
		The following is a short exact sequence:
		\begin{equation}
			0 \longrightarrow \mathfrak{m}_\K \xrightarrow{\Expx{\edwards}} \qpjacobianzof{\K}{\edwards} \ \xrightarrow{\Modx{\edwards}} \qpjacobianzof{k}{\edwards} \longrightarrow 0.
		\end{equation}
	\end{corollary}
	\begin{proof}
		The proof follows from the fact that, from \thmref{exp-edwards}, $ \Expx{\edwards} $ is a monomorphism
		and $ \Modx{\edwards} $ is an epimorphism. Moreover,
		since, for any $ z \in \mathfrak{m}_\K $, $ \Expxof{\edwards}{z} = (\mathcal{O}(z^3), 1 + \mathcal{O}(z^3)) $,
		then $ \Modxof{\edwards}{\Expxof{\edwards}{z}} = (0, 1) $,
		and $ \Imgof{\Expx{\edwards}} \subseteq \Kerof{\Modx{\edwards}} $. Finally,
		together with $ \Expx{\redweier} $, which is invertible by \eqref{eq:inv-exp},
		the map $ \Expx{\edwards} = \beta \circ \chi \circ \Expx{\redweier} $ is invertible
		for $ z \in \mathfrak{m}_\K $,
		that is, one can write any point $ P $ in $ \Kerof{\Modx{\edwards}} $ as $ P = \Expxof{\edwards}{z} $,
		for some $ z \in \mathfrak{m}_\K $, thus
		$ \Kerof{\Modx{\edwards}} \subseteq \Imgof{\Expx{\edwards}} $.
	\end{proof}

	In the following, we stress the meaning of choosing a curve whose cardinality differs from
	the cardinality of its ground field.

	\begin{definition}
		Let $ \weier $ be an elliptic curve such that $ \cardof{\qpjacobianof{k}{\weier}} = \cardof{k} $,
		where $ k $ is a finite field. The curve $ \weier $ is an \keyphrase{anomalous curve}.
	\end{definition}

	Non-anomalous curves are subject to attacks by means, for instance, of pairing mappings,
	that is, efficiently computable, bilinear, and non-degenerate
	maps $ e \colon G_1 \times G_2 \rightarrow G_3 $, where
	typically $ G_1 $ and $ G_2 $ are cyclic subgroups (such as the Weil pairing,
	see e.g. $ \S $III.8 in \cite{Silverman2009}, and the Ate pairing, see \cite{AtePairing})
	or quotient groups (such as
	the Tate pairing, see e.g. \cite{FreyMullerRuck}, and the Eta pairing, see \cite{EtaPairing})
	of the Jacobian of the curve, while $ G_3 $ is a subgroup of the multiplicative group
	of the ground field because the pairing carries the logarithm of an element in $ G_1 $
	to the logarithm of an element in $ G_3 $
	(see e.g. \cite{MenezesOkamotoVanstone}).
	Anomalous curves are safe with respect to these attacks since
	all the above pairings are defined if and only if the cardinalities
	of $ G_1 $ and $ G_2 $ divide $ q^k - 1 $,
	where $ q^k $ is the cardinality of the ground field.

	On the other hand, however, anomalous curves are also subject to attacks
	as it is possible to map the Jacobian of such curves to
	the additive group of the finite field $ k $
	(see \cite{LEPREVOST2005225, MR1863606, Smart1999, Silverman2009, SamaevAnomalous, SatohArakiAnomalous}).
	We address the reader to $ \S $XI.6 in \cite{Silverman2009} for a simple polynomial algorithm able to solve the
	ECDLP for an anomalous curve.

	\begin{theorem}\label{thm:jw-splitting}
		If $ k = \mathcal{O}_\K / \mathfrak{m}_\K $ is finite, and $ \weier $ is not an anomalous curve,
		then $ \qpjacobianof{\K}{\weier} $ is isomorphic to the direct sum of
		$ \qpjacobianof{k}{\weier} $ and $ \mathfrak{m}_\K $.
	\end{theorem}

	\begin{proof}
		As $ k = \mathcal{O}_\K / \mathfrak{m}_\K $ is finite and $ \weier $ is not anomalous,
		for any $ 1 \le h \in \Z $,	the sequence:
		\begin{equation*}
			0 \longrightarrow H \longrightarrow \qpjacobianof{H}{\weier} \longrightarrow \qpjacobianof{k}{\weier} \longrightarrow 0,
		\end{equation*}
		where $ H = \mathfrak{m}_\K / (\varpi^h_\K \mathcal{O}_\K) $,
		with $ \varpi_\K $ the uniformizer of $ \K $, is splitting by the Schur-Zassenhaus theorem
		and defines, therefore, a section $ \sigma^h_{\weier} \colon \qpjacobianof{k}{\weier} \rightarrow \qpjacobianof{H}{\weier} $ which is a homomorphism. Taking the inverse limit
		$ \sigma_{\weier} = \lim\limits_{h \rightarrow \infty} \sigma^{h}_{\weier} $,
		we obtain a section $ \sigma_{\weier} \colon \qpjacobianof{k}{\weier} \rightarrow \qpjacobianof{H}{\weier} $
		which is a homomorphism, hence the sequence is splitting.
	\end{proof}

	\begin{corollary}\label{corollary:exp-split-exact-sequence}
		If $ k = \mathcal{O}_\K / \mathfrak{m}_\K $ is finite, $ \weier $ is not an anomalous curve,
		and $ \edwards $ is the Edwards curve birational equivalent to $ \weier $,
		then $ \qpjacobianzof{\K}{\edwards} $ is isomorphic to $ \qpjacobianzof{k}{\edwards} \oplus \mathfrak{m}_\K $.
	\end{corollary}
	\begin{proof}
		Since we confined ourselves to the case in \thmref{isomorphism-w-e}, then
		we have that
		$ \qpjacobianzof{\K}{\edwards} \cong \qpjacobianof{\K}{\weier} $,
		$ \qpjacobianzof{k}{\edwards} \cong \qpjacobianof{k}{\weier} $,
		and the proof follows from \thmref{jw-splitting}.
	\end{proof}

	Note that the above exact sequence does not split over $ \K $
	if one supposes that the elliptic curve in Weierstrass form,
	taken into account in \thmref{isomorphism-w-e}, is an \keyphrase{anomalous curve},
	as we show in the next example.

	\begin{example}
		Let $ \weier $ be the elliptic curve in Weierstrass form
		defined by the equation $ y^2 = x^3 + 4 x + 7 $ over $ k = \GF{53} $,
		whose Jacobian can be readily verified to have $ 53 $ elements. Hence,
		$ \qpjacobianof{k}{\weier} $ is isomorphic to the cyclic group $ C_{53} $.
		However, $ \qpjacobianof{\Z / {53}^2 \Z }{\weier} \ne C_{53} \oplus C_{53} $
		as the point $ P = (3, 130) \in \weierof{\Zof{{53}^2}} $
		is such that $ 53 (P - \Omega) = \homcoords{0}{53}{1603} - \Omega \ne \Omega - \Omega $.
	\end{example}

	\section{The map \texorpdfstring{$ \Exp $}{Exp} for Edwards curves over \texorpdfstring{$ \Qp{p} $}{Qp}}\label{sec:3}

	The goal of this section is to compute the map $ \Exp $ for Edwards curves $ \edwards $
	over the local field $ \Qp{p} $ of $ p $-adic numbers.

	In particular, we study the field $ \Qp{p} $ through the inverse limit
	$ \Z_p = \lim\limits_{\longleftarrow} \Zof{p^k} $, that is,
	we compute $ \Zof{p} $, $ \Zof{p^2} $, $ \ldots $, $ \Zof{p^k} $
	approaching $ \Z_p $, the field of $ p $-adic integers, for $ k \rightarrow \infty $.

	The field $ \Qp{p} $ is a non-Archimedean local field of characteristic zero, thus
	putting $ \K = \Qp{p} $, one has that its ring of integers $ \mathcal{O}_\K $
	is the ring $ \Z_p $ of $ p $-adic integers, its prime ideal $ \mathfrak{m}_\K $
	is $ p \Z_p $ (which uniformizer $ \varpi_\K $ is equal to $ p $), and its residue field
	$ \mathcal{O}_\K / \mathfrak{m}_\K $ is $ \Zof{p} = \GF{p} $.

	Recall that, in \remarkref{series-expansion}, we gave the Laurent
	series expansion for the Weierstrass $ \wp $-function and its
	derivative $ \wp^\prime $ for a complex number $ z $.
	As here we are now focusing on the field $ \Qp{p} $,
	one has to take into account the convergence radius of these series over $ \Qp{p} $.
	In the context of the field of $ p $-adic numbers,
	a convergence neighborhood of zero
	is given by multiples of $ p $, that is, when $ p \mid z $.
	In this neighborhood, these series always converge since $ c_k z^{2k - 2} \equiv 0\ \modp{p^h} $,
	and $ (2k - 2) c_k z^{2k - 3} \equiv 0\ \modp{p^h} $, for a suitable positive integer $ h $.

	If $ \qpjacobianzof{k}{\edwards} $ denotes the image of the subgroup $ \jacobianzof{\edwards} $
	modulo $ p^k $, then, applying the above changes and observations,
	the results in \secref{2} can be expressed as follows:
	$ \qpjacobianzof{k}{\edwards} = \qpjacobianzof{1}{\edwards} \oplus \Imgof{\Expx{\edwards}} $.

	Indeed, the map $ \Expx{\edwards} $ over $ \Qp{p} $ can be expressed, through an inverse limit
	with $ k $ increasing, as follows:
	\begin{alignat*}{2}
		\Expx{\edwards} \colon \  &\frac{p \Z_p}{p^k \Z_p} &\ \longrightarrow\  &\edwardsof{\Zof{p^k}}\\
			&z = ph &\ \longmapsto\  & \roundbrack{\frac{2}{3}\frac{y_1 (3 \wp(z) - a^\prime)}{x_1 \wp^\prime(z)}, \frac{3 \wp(z) - a^\prime - 3 x_1}{3 \wp(z) - a^\prime + 3 x_1}},
	\end{alignat*}
	where $ h = 1, 2, \ldots, p^{k - 1} $, since $ \Z_p = \lim\limits_{\longleftarrow} \Zof{p^k} $.

	\begin{remark}
		Note that in this case we consider, as the domain of $ \Expx{\edwards} $,
		the quotient $ \frac{p \Z_p}{p^k \Z_p} $ since, modulo $ p^k $,
		$ \Expxof{\edwards}{ph} = \Expxof{\edwards}{p (h + p ^ k)}$, with $ h \in \Z $.
	\end{remark}

	\begin{remark}
		Note that, as there is a natural isomorphism from $ \Imgof{\Exp} $ and $ \Zof{p^{k - 1}} $
		through the following map:
		\begin{alignat*}{2}
			\frac{p \Z_p}{p^k \Z_p} &\ \longrightarrow\ &\ \Zof{p^{k - 1}}\\
			ph &\ \longmapsto\  &\ h,
		\end{alignat*}
		where $ h = 1, 2, \ldots, p^{k - 1} $,
		we have that $ \qpjacobianzof{k}{\edwards} = \qpjacobianzof{1}{\edwards} \oplus \Zof{p^{k - 1}} $.
	\end{remark}

	It is necessary to make some further adjustments to the map $ \Exp $ for elliptic curves in Weierstrass
	form. In particular, as we are approximating $ \Qp{p} $ with $ \Zof{p^k} $, with $ k \rightarrow \infty $,
	the map $ \Exp $ for elliptic curves in short Weierstrass form should be rewritten as follows:
	\begin{alignat*}{2}
		\Expx{\redweier} \colon \  &\frac{p \Z_p}{p^k \Z_p} &\ \longrightarrow\  &\weierof{\Zof{p^k}}\\
		&z &\ \longmapsto\  & \homcoords{t z^3}{t z^3 \wp(z)}{t z^3\frac{1}{2} {\wp^\prime(z)}},
	\end{alignat*}
	where $ t $ is the least common multiple between the denominators of the series expansion
	of $ \wp $ and $ \wp^\prime $ (see \remarkref{series-expansion}). In particular,
	the multiplication by the factor $ t z^3 $ has to be done in order to make all coordinates integer,
	and therefore to avoid modular inversions when the denominator is a multiple of $ p $
	as we move from the field $ \Zof{p} $ to the ring $ \Zof{p^k} $.

	For elliptic curves in short Weierstrass form,
	since $ \Imgof{\Expx{\redweier}} = \Kerof{\Modx{\redweier}} $,
	then $ P \in \Imgof{\Expx{\redweier}} $ if $ \Modxof{\redweier}{P} = \Omega $.
	Therefore, the points belonging to $ \Imgof{\Expx{\redweier}} $ have the following form
	$ P = \homcoords{p h_1}{p h_2}{h_3} $ with $ p \nmid h_3 $.
	However, since $ \Omega $, the point at infinity of $ \weier $, is mapped through $ \beta $
	onto the neutral point $ O \in \edwardsof{\Zof{p^k}} $,
	then all the points belonging to $ \Imgof{\Expx{\edwards}} $ will be equivalent,
	modulo $ p $, to $ O $ and, therefore, they will all be affine points.
	So, by counting all the affine points $ (x, y) \in \edwardsof{\Zof{p^k}} $,
	one may check that the number of these points is equal to $ \cardof{\qpjacobianzof{1}{\edwards}} \cdot p^{k - 1} $,
	where $ p^{k - 1} $ is the cardinality of $ \Imgof{\Expx{\edwards}} $.

	Thus, the map $ \Exp $ allows us to speed up the addition operation by splitting
	the original group $ \qpjacobianzof{k}{\edwards} $ into a pair $ (P, c) $, where $ P \in \edwardsof{\Zof{p}} $
	and $ c \in \Zof{p^{k - 1}} $.

	\begin{remark}
		Note that the addition formula for the Weierstrass form cannot be applied in
		the case in which points, reduced modulo $ p $, return the point at infinity.
		In this latter case, the sum of two points over $ \Zof{p^k} $ may return a point that does not exist,
		such as $ \homcoords{0}{0}{0} $.
		On the contrary, every point in $ \edwardsof{\Zof{p^k}} $ is affine, thus the addition formula
		for an Edwards curve always returns the proper result.
	\end{remark}

	\section{Conclusions}

	The Edwards curves are a recent ($ 2007 $) mathematical tool used in cryptographic
	and digital signature applications because of their efficient (and secure) group law operations.
	Until now, these curves have been studied in order to find other ways to employ them and
	further speed up their applications, while preserving their security.

	In this paper, we extended the map $ \Exp $ for elliptic curves in short
	Weierstrass form, defined over $ \C $ by the equation $ y^2 = x^3 + a x + b $, to
	the Edwards curves $ \edwards $, defined over local fields by the equation
	$ {\hat{x}}^2 + {\hat{y}}^2 = 1 + d {\hat{x}}^2 {\hat{y}}^2 $, with $ d $ a non-square,
	by using the birational equivalence between the Weierstrass form $ \weier $
	of equation $ y^2 = x^3 + a^\prime x^2 + b^\prime x $ and $ \edwards $.

	Up to the representation of the elements of $ \qpjacobianzof{\K}{\edwards} $
	as pairs $ (P, c) $, where $ P \in \edwardsof{\GF{p}} $
	and $ c \in \Zof{p^{k - 1}} $, this map provides a tool able to speed up the
	group law operations for Edwards curves over local fields, and in particular over the
	field $ \Qp{p} $ of $ p $-adic numbers
	by splitting the whole subgroup $ \qpjacobianzof{k}{\edwards} $,
	that is, the image of the subgroup $ \jacobianzof{\edwards} $ over $ \Zof{p^k} $,
	into a pair $ (P, c) $, where $ P \in \edwardsof{\GF{p}} $
	and $ c \in \Zof{p^{k - 1}} $.

	This also gives a motivation for studying the map able to correctly define the pair
	$ (P, c) $ and a group law in order to sum two of such elements, that is,
	$ (P_1, c_1) + (P_2, c_2) $.

	\printbibliography[title={References}]

@misc{GFGoppaCodes,
	author={Giuseppe Filippone},
	title={Goppa codes over Edwards curves},
	day={},
	month={},
	year={},
	note={(in progress)},
	url={https://drive.google.com/file/d/120dQ2xZ-Oz0XmCBQn3rSvIGMg6K4dBmK/view?usp=share_link}
}

@article{kosterspannekoek,
	author = {Kosters, Michiel and Pannekoek, René},
	year = {2017},
	month = {03},
	pages = {},
	title = {On the structure of elliptic curves over finite extensions of $ \Qp{p} $ with additive reduction},
	journal = {arXiv:1703.07888},
	doi = {https://doi.org/10.48550/arXiv.1703.07888}
}

@article{Tang2015,
	author={Tang, ChunMing
	and Xu, MaoZhi
	and Qi, YanFeng},
	title={Cryptography on twisted Edwards curves over local fields},
	journal={Science China Information Sciences},
	year={2015},
	month={01},
	day={01},
	volume={58},
	number={1},
	pages={1-15},
	issn={1869-1919},
	doi={10.1007/s11432-014-5155-z},
	url={https://doi.org/10.1007/s11432-014-5155-z}
}

@article{Xu2008,
	author={Xu, MaoZhi
	and Zhao, ChunLai
	and Feng, Min
	and Ren, ZhaoRong
	and Ye, JiQing},
	title={Cryptography on elliptic curves over p-adic number fields},
	journal={Science in China Series F: Information Sciences},
	year={2008},
	month={03},
	day={01},
	volume={51},
	number={3},
	pages={258-272},
	issn={1862-2836},
	doi={10.1007/s11432-008-0014-4},
	url={https://doi.org/10.1007/s11432-008-0014-4}
}

@article{XuYue,
	author={Yue, Zhi Hong and Xu, Mao Zhi},
	title={Hierarchical Management Scheme by Local Fields},
	journal={Acta Mathematica Sinica},
	year={2010},
	month={12},
	day={15},
	volume={27},
	number={1},
	pages={155-168},
	issn={1439-8516},
	url={https://actamath.cjoe.ac.cn/Jwk_sxxb_en/EN/10.1007/s10114-011-9110-2},
	doi={https://doi.org/10.1007/s10114-011-9110-2}
}

@article{AtePairing,
	author={Hess, F. and Smart, N.P. and Vercauteren, F.},
	journal={IEEE Transactions on Information Theory},
	title={The Eta Pairing Revisited},
	year={2006},
	volume={52},
	number={10},
	pages={4595-4602},
	doi={10.1109/TIT.2006.881709}
}

@article{EtaPairing,
	author={Barreto, Paulo S. L. M.
	and Galbraith, Steven D.
	and h{\'E}igeartaigh, Colm {\'O}'
	and Scott, Michael},
	title={Efficient pairing computation on supersingular Abelian varieties},
	journal={Designs, Codes and Cryptography},
	year={2007},
	month={03},
	day={01},
	volume={42},
	number={3},
	pages={239-271},
	issn={1573-7586},
	doi={10.1007/s10623-006-9033-6},
	url={https://doi.org/10.1007/s10623-006-9033-6}
}

@book{SerreLocalField,
	AUTHOR = {Serre, Jean-Pierre},
	TITLE = {Local fields},
	SERIES = {Graduate Texts in Mathematics},
	VOLUME = {67},
	NOTE = {Translated from the French by Marvin Jay Greenberg},
	PUBLISHER = {Springer-Verlag},
	ADDRESS = {New York},
	YEAR = {1979},
	PAGES = {viii+241},
	ISBN = {0-387-90424-7},
	MRCLASS = {12Bxx},
	MRNUMBER = {554237 (82e:12016)},
	BOEKCODE = {11Sxx},
	DOI = {https://doi.org/10.1007/978-1-4757-5673-9}
}

@article{SamaevAnomalous,
	author = {Semaev, I. A.},
	title = {Evaluation of Discrete Logarithms in a Group of P-Torsion Points of an Elliptic Curve in Characteristic p},
	year = {1998},
	issue_date = {Jan. 1998},
	publisher = {American Mathematical Society},
	address = {USA},
	volume = {67},
	number = {221},
	issn = {0025-5718},
	url = {https://doi.org/10.1090/S0025-5718-98-00887-4},
	doi = {10.1090/S0025-5718-98-00887-4},
	journal = {Math. Comput.},
	month = {01},
	pages = {353–356},
	numpages = {4}
}

@inproceedings{SilvermanLiftingECDLP,
	author="Silverman, Joseph H.",
	editor="Avanzi, Roberto Maria
	and Keliher, Liam
	and Sica, Francesco",
	title="Lifting and Elliptic Curve Discrete Logarithms",
	booktitle="Selected Areas in Cryptography",
	year="2009",
	publisher="Springer Berlin Heidelberg",
	address="Berlin, Heidelberg",
	pages="82--102",
	isbn="978-3-642-04159-4",
	url={https://doi.org/10.1007/978-3-642-04159-4_6},
	doi={10.1007/978-3-642-04159-4_6}
}

@misc{SatohArakiAnomalous,
	author = {Satoh, Takakazu and Araki, Kiyomichi},
	title = {Fermat Quotients and the Polynomial Time Discrete Log Algorithm for Anomalous Elliptic Curves},
	year = {1998},
	volume = {47},
	url = {https://doi.org/10.14992/00009878},
	doi = {10.14992/00009878},
	journal = {Commentarii mathematici Universitatis Sancti Pauli},
	month = {06},
	pages = {81-92}
}

@article {MR1863606,
	AUTHOR = {Jacobson, Michael and Menezes, Alfred and Stein, Andreas},
	TITLE = {Solving elliptic curve discrete logarithm problems using {W}eil descent},
	JOURNAL = {J. Ramanujan Math. Soc.},
	FJOURNAL = {Journal of the Ramanujan Mathematical Society},
	VOLUME = {16},
	YEAR = {2001},
	NUMBER = {3},
	PAGES = {231--260},
	ISSN = {0970-1249},
	MRCLASS = {14G50 (11G20 11Y16)},
	MRNUMBER = {1863606},
	MRREVIEWER = {Steven D. Galbraith},
	URL = {https://eprint.iacr.org/2001/041}
}

@ARTICLE{FreyMullerRuck,
	author={Frey, G. and Muller, M. and Ruck, H.-G.},
	journal={IEEE Transactions on Information Theory},
	title={The Tate pairing and the discrete logarithm applied to elliptic curve cryptosystems},
	year={1999},
	volume={45},
	number={5},
	pages={1717-1719},
	doi={10.1109/18.771254}
}

@inbook{Lange2011,
	author="Lange, Tanja",
	editor="van Tilborg, Henk C. A.	and Jajodia, Sushil",
	title="Edwards Curves",
	bookTitle="Encyclopedia of Cryptography and Security",
	year="2011",
	publisher="Springer US",
	address="Boston, MA",
	pages="380--382",
	isbn="978-1-4419-5906-5",
	doi="10.1007/978-1-4419-5906-5_243"
}

@incollection{AbramowitzStegun1972,
	author = {Abramowitz, Milton and Stegun, Irene A.},
	title = {Weierstrass Elliptic and Related Functions},
	booktitle = {Handbook of Mathematical Functions with Formulas, Graphs and Mathematical Tables},
	editor = {Abramowitz, Milton and Stegun, Irene A.},
	chapter = {18},
	pages = {627-671},
	address = {New York},
	publisher = {Dover Publications, Inc.},
	year = {1972},
	isbn = {978-0486612720}
}

@article{ARENE2011842,
	title = {Faster computation of the Tate pairing},
	journal = {Journal of Number Theory},
	volume = {131},
	number = {5},
	pages = {842-857},
	year = {2011},
	note = {Elliptic Curve Cryptography},
	issn = {0022-314X},
	doi = {https://doi.org/10.1016/j.jnt.2010.05.013},
	url = {https://www.sciencedirect.com/science/article/pii/S0022314X10001757},
	author = {Christophe Arène and Tanja Lange and Michael Naehrig and Christophe Ritzenthaler}
}

@incollection{ApostolTomMike1990MfaD,
	series = {Graduate texts in mathematics 41},
	publisher = {Springer Verlang},
	booktitle = {Modular functions and Dirichlet series in number theory},
	year = {1996},
	edition = {2},
	address = {New York},
	author = {Apostol, Tom Mike},
	isbn = {978-0387971278},
	chapter = {1.6 - 1.11},
	pages = {9 - 14}
}

@article{MenezesOkamotoVanstone,
	author={Menezes, A.J. and Okamoto, T. and Vanstone, S.A.},
	journal={IEEE Transactions on Information Theory},
	title={Reducing elliptic curve logarithms to logarithms in a finite field},
	year={1993},
	volume={39},
	number={5},
	pages={1639-1646},
	doi={10.1109/18.259647}
}

@article{Smart1999,
	author={Smart, N. P.},
	title={The Discrete Logarithm Problem on Elliptic Curves of Trace One},
	journal={Journal of Cryptology},
	year={1999},
	month={06},
	day={01},
	volume={12},
	number={3},
	pages={193-196},
	issn={1432-1378},
	doi={10.1007/s001459900052},
	url={https://doi.org/10.1007/s001459900052}
}

@article{LEPREVOST2005225,
	title = {Generating anomalous elliptic curves},
	journal = {Information Processing Letters},
	volume = {93},
	number = {5},
	pages = {225-230},
	year = {2005},
	issn = {0020-0190},
	doi = {https://doi.org/10.1016/j.ipl.2004.11.008},
	url = {https://www.sciencedirect.com/science/article/pii/S0020019004003527},
	author = {Franck Leprévost and Jean Monnerat and Sébastien Varrette and Serge Vaudenay}
}

@book{Silverman2009,
	series = {Graduate Texts in Mathematics},
	volume = {106},
	publisher = {Springer Verlang},
	title = {The Arithmetic of Elliptic Curves},
	year = {2009},
	edition = {2},
	address = {New York},
	author = {Silverman, Joseph H.},
	isbn = {978-0-387-09493-9},
	doi = {10.1007/978-0-387-09494-6}
}

@inproceedings{HisilWongCarterDawson,
	author = {Hisil, Huseyin and Wong, Kenneth Koon-Ho and Carter, Gary and Dawson, Ed},
	title = {Faster Group Operations on Elliptic Curves},
	year = {2009},
	isbn = {9781920682798},
	publisher = {Australian Computer Society, Inc.},
	address = {AUS},
	doi = {10.5555/1862758.1862762},
	url = {https://dl.acm.org/doi/10.5555/1862758.1862762},
	booktitle = {Proceedings of the Seventh Australasian Conference on Information Security - Volume 98},
	pages = {7–20},
	numpages = {14},
	location = {Wellington, New Zealand},
	series = {AISC '09}
}

@article{Edwards,
	author = {Edwards, Harold},
	year = {2007},
	month = {07},
	pages = {393-423},
	title = {A normal form for elliptic curves},
	volume = {44},
	journal = {Bulletin of The American Mathematical Society - BULL AMER MATH SOC},
	doi = {10.1090/S0273-0979-07-01153-6}
}

@inproceedings{BernsteinBirknerLangePeters,
	author="Bernstein, Daniel J. and Birkner, Peter and Lange, Tanja and Peters, Christiane",
	editor="Srinathan, K. and Rangan, C. Pandu and Yung, Moti",
	title="Optimizing Double-Base Elliptic-Curve Single-Scalar Multiplication",
	booktitle="Progress in Cryptology -- INDOCRYPT 2007",
	year="2007",
	publisher="Springer Berlin Heidelberg",
	address="Berlin, Heidelberg",
	pages="167--182",
	isbn="978-3-540-77026-8",
	doi="10.1007/978-3-540-77026-8_13"
}

@inproceedings{BernsteinLange2,
	author="Bernstein, Daniel J. and Lange, Tanja",
	editor="Bozta{\c{s}}, Serdar
	and Lu, Hsiao-Feng (Francis)",
	title="Inverted Edwards Coordinates",
	booktitle="Applied Algebra, Algebraic Algorithms and Error-Correcting Codes",
	year="2007",
	publisher="Springer Berlin Heidelberg",
	address="Berlin, Heidelberg",
	pages="20--27",
	isbn="978-3-540-77224-8",
	doi="10.1007/978-3-540-77224-8_4"
}

@inproceedings{BernsteinLange,
	author="Bernstein, Daniel J. and Lange, Tanja",
	editor="Kurosawa, Kaoru",
	title="Faster Addition and Doubling on Elliptic Curves",
	booktitle="Advances in Cryptology -- ASIACRYPT 2007",
	year="2007",
	publisher="Springer Berlin Heidelberg",
	address="Berlin, Heidelberg",
	pages="29--50",
	isbn="978-3-540-76900-2",
	doi="10.1007/978-3-540-76900-2_3"
}

\end{document}